\documentclass[12pt]{article}
\usepackage{amsmath,amsfonts,amssymb,amsthm}

\def\noi{\noindent}

\def\pf{\noi{\bf Proof.\ \,}}
\def\qed{{$\square$}}

\def\em{\it }        

\def\l{\lambda}
\def\o{\omega}

\def\C{{\mathbb C}}

\def\Z{{\mathbb Z}}

\def\tr{{\rm tr}}

\def\la{\langle}
\def\ra{\rangle}
\def\<{\langle}
\def\>{\rangle}

\def\bs{\it}            
\def\Aut{{\bs Aut}}
\def\dim{{\bs dim}}

\def\l{{\lambda}}
\def\Res{{\rm Res}}
\def\wt{{\rm wt}}

\def\vac{\hbox{\bf 1}} 

\begin{document}
\begin{center}
{\Large 
Automorphism groups and derivation algebras of finitely 
generated vertex operator algebras}

\vspace{10mm}
Chongying Dong\footnote{The first author is supported by NSF grant 
DMS-9987656 and a research grant from the Committee on Research, UC Santa Cruz.}
\\[0pt]
Department of Mathematics, University of California\\[0pt]
Santa Cruz, CA 95064\\[0pt]
\vspace{3mm}
Robert L.~Griess Jr.\footnote{The second author is supported by 
NSA grant USDOD-MDA904-00-1-0011.} \\[0pt]
Department of Mathematics, University of Michigan\\[0pt]
Ann Arbor, MI 48109
\end{center}

\newtheorem{thm}{Theorem}[section]
\newtheorem{prop}[thm]{Proposition}
\newtheorem{lem}[thm]{Lemma}
\newtheorem{rem}[thm]{Remark}
\newtheorem{coro}[thm]{Corollary}
\newtheorem{conj}[thm]{Conjecture}
\newtheorem{de}[thm]{Definition}

\newtheorem{nota}[thm]{Notation}
\newtheorem{ex}[thm]{Example}


\begin{abstract}

We investigate the general structure of the automorphism group and 
the Lie algebra of
derivations of a finitely generated vertex operator algebra.   
The automorphism group is isomorphic to an algebraic group.
Under natural assumptions, the derivation algebra
has an invariant bilinear form and the 
ideal of inner derivations is nonsingular. 
\end{abstract}

\section{Introduction}

A study of the automorphism group of a vertex operator algebra
(VOA)  is a natural part of the study of 
vertex operator algebras. 

\begin{de}{\rm 
Let $V$ be  a vertex operator algebra.  
We say that $a \in GL(V)$ is an {\it automorphism of $V$} if and only if it 
leaves the vacuum element and the principal Virasoro element fixed
($a\vac = \vac$ and $a\omega = \omega$) and 
preserves all $V$-compositions, i.e., for all 
$m \in \Z$ and $u, v \in V$, we have $a(u_mv)=a(u)_ma(v)$.  
It follows that an automorphism fixes all the 
$V_i$ since they are eigenspaces for an operator in the 
series for the principal Virasoro element.  

The set of all automorphisms is a group, denoted
$Aut(V)$.
}
\end{de}

In the definition, it suffices to restrict $u$
and
$v$  to homogeneous elements. 

In some definitions of VOA automorphism, there is no requirement
that the principal Virasoro element be fixed.  

So far, we know the automorphism groups explicitly 
for relatively few  vertex operator
algebras, such as $V^{\natural}$ [FLM], vertex operator algebra
$V_L$ for a  positive definite even lattice 
$L$ [DN], certain vertex operator algebras with central charge 1 
([DG], [DGR]), vertex operator algebras associated to highest weight 
representations for affine algebras (cf. [DLY]), vertex 
operator algebras associated to codes [M] and a few special 
cases, e.g., [G].

The determination
of each of these automorphism groups has  its own story and depends
heavily on the specifics of  the auxiliary object used to construct the
VOA,  such as a
lattice, Lie algebra, or code. Nevertheless, one can observe 
that all these automorphism
groups have similarities. 

We denote by 
$(V, k^{th})$ the algebra with underlying vector space 
$V$ and product $a_kb$, for $a, b \in V$, where $a_k$ is the
coefficient at 
$z^{-k-1}$ in the vertex operator for $a$.  The linear subspace $V_m$ is
closed under this product if $m-k=1$.  This algebra is denoted 
$(V_m, (m-1)^{th})$.  The cases $m=1$ and 2 are especially interesting.

For simplicity of exposition, let us assume that $V$ is a simple vertex
operator algebra of CFT type (see [DLMM], or definition 2.8 of 
this article). 
Then $(V_1,0^{th})$ is a Lie algebra with bracket $[u,v]=u_0v$
for $u,v\in V_1$ where $Y(u,z)=\sum_{n\in\Z}u_nz^{-n-1}.$  
The endomorphism 
$u_0$ is a derivation of $V$ in the sense
that $u_0{\bf 1}=0,$ $u_0\omega=0$ and 
$$u_0(Y(v,z)w)=Y(u_0v,z)w+Y(u,z)u_0w$$ for any $v,w\in V$
(cf. [DN]). Moreover,
the exponential $e^{u_0}$ is an automorphism of $V$ ($e^{u_0}$ is a
well defined operator on $V$ as each $V_n$ is finite dimensional).  
Denote by $Aut_1(V)$ 
the subgroup of the automorphism group $Aut(V)$ of $V$ generated by
$e^{u_0}$ for $u\in V_1.$ Then $Aut_1(V)$ is a finite dimensional
connected algebraic normal subgroup of $\Aut (V)$.  
In all the examples discussed above, $\Aut(V)/Aut_1(V)$ is a finite
group.  We think that this is probably a general phenomenon for 
rational vertex operator algebras.  There are counterexamples
when $V$ is not rational (cf. [DM1] and Examples \ref{exama} and 4.1 in this
article). 

It has been a feeling 
 for a while that any rational vertex operator 
algebra of CFT type is finitely generated. If $V$ is regular  in the sense 
that any weak module is a direct sum of
ordinary modules (see [DLM2]),  it is proved in [KL] and [L3] (also see [GN]) 
that $V$ is finitely generated. 
It is also felt that rational vertex operator algebras must be  
regular.  Interest in the category of modules 
is motivation 
to study automorphism groups of
finitely generated vertex operator algebras. 

The main result of this paper is that the automorphism group $Aut(V)$ 
of a finitely generated vertex operator algebra is isomorphic
to a finite
dimensional algebraic group. It is well known that a finite dimensional
algebraic group $G$ has only finitely
many connected components and so $G/G^0$ is a finite
group, where $G^0$ is the connected component of $G$ containing the 
identity. We expect that the normal subgroup
$Aut_1(V)$ of $Aut(V)$ is exactly $Aut(V)^0$ for all rational vertex
operator algebras $V$ of CFT type. This 
property holds for all examples discussed before.

There is a close relation between the automorphism group and the
Lie algebra of derivations of a vertex operator algebra. If $d$ 
is a derivation of a vertex operator algebra $V$ then $e^d$
is an automorphism of $V$ (see Section 3).  If $V$ is a finitely generated rational vertex operator algebra of CFT type,
equality of $Aut_1(V)$ and $Aut(V)^0$ is equivalent to 
all the derivations of $V$ being given by $u_0$ for $u\in V_1.$

The paper is organized as follows. In Section 2 we prove that
the automorphism group of a finitely generated vertex operator algebra 
is a finite dimensional algebraic group. We also give an example 
of a non finitely generated vertex operator algebra
whose automorphism group is not isomorphic to 
an algebraic group. In Section
3 we study derivations of vertex operator algebras. For
$v\in V$ we define a linear operator $o(v)$ by the conditions that 
$o(v)=v_{\wt v-1}$ if $v$ is homogeneous.  
We show in Section 3 that $o(v)$ is a derivation of $V$ if and only
if $v\in V_1.$  We also show that the Lie algebra $V_1$ is
an ideal of the Lie algebra of the derivations and has an orthogonal
complement with respect to a suitable invariant symmetric bilinear form.
In Section 4 we discuss
an example of a nonsimple finitely generated 
vertex operator algebra and its automorphism
group. 

\section{Automorphism groups}

We suppose that the VOA $V$ is finitely generated (cf. [FHL]). 
This is equivalent to existence
of $n\in\Z$ so that $U=\oplus_{m\leq n}V_m$ generates $V$ in the sense that
$$V=span\<u^1_{i_1}\cdots u^s_{i_s}u \ | 
\ u^j,u\in U, s \in \{0,1,2,\dots \}, i_j\in\Z\>.$$

For a subset $A$ of $V$, set $A^{r+1} := A \times \cdots \times A$
($r+1$ times) and $A^\infty := \bigcup_{r \ge 0}A^{r+1}$. An element
of $A^\infty$ is a finite length vector $\vec x = (x_0, \dots , x_r)$
and we call $r+1$ the {\it length} of $\vec x$. For every nonempty
finite sequence $\vec m := (m_1, \dots , m_r)$ of integers, we define
the function, called $\vec m$-{\it composition} $\mu :=\mu_{\vec m} :
V^{r+1} \rightarrow V$ by $\mu (x_0, \dots ,x_r) :=
(x_0)_{m_1}(x_2)_{m_2} \cdots (x_{r-1})_{m_{r}}x_r$. Call $r+1$ the
{\it length} of $\mu$.  Thus, $\mu (\vec x)$ is defined if and only if
$\mu$ and $\vec x$ have the same lengths, in which case we say that $(\mu,
\vec x)$ is an {\it admissible pair}. In case the entries of such
$\vec x$ are restricted to a subset $A$ of $V$, we all the pair an
{\it $A$-admissible pair}.  If the coordinates of $\vec x = (x_0,
\dots , x_r)$ are homogeneous, then $\mu (\vec x)$ is homogeneous and
we define the {\it weight} of $(\mu, \vec x)$ to be $\sum_{i=0}^r
wt(x_i) - \sum_{i=1}^r (m_i + 1)$. This is just the weight of $\mu
(\vec x)$ if $\mu (\vec x)$ is nonzero.

Such a function, for some $\vec m$, is called a $V$-{\it composition}. 
For a subset $A$ of $V$, the restriction
of $\mu$ to tuples of elements in $A$ is denoted $\mu_A$.

\begin{rem} {\rm The property that $U$ generates $V$ means that for each integer $m$ there exists
a finite set $S$ of $V$-compositions so that $V_m  = \sum_{\mu \in S} Im(\mu_U) \cap V_m$.}
\end{rem}

We choose a basis $\cal B$ of
$U$ consisting of homogeneous elements, including $\vac$.
Let $\cal Q$ be the set of $\cal B$-admissible pairs.  Define ${ \cal
Q}_m$ to be the set of pairs in $\cal Q$ of weight $m$.  Then $V_m$ is
spanned by a finite set of ``monomials"  in $\cal B$ of weight $m$, i.e.,
elements of certain $Im(\mu_{\cal B})$.

There is a finite set ${\cal R}_m$ of pairs $(\mu, \vec x) \in {\cal Q}_m$ so that the
set ${\cal B}_m := \{ \mu (\vec x) \ | \ (\mu, \vec x) \in {\cal R}_m \}$ forms a basis
for $V_m$. Choose ${\cal R}_0=\{ (\mu_0, \vac ) \}$ ($\mu_0$ is the
trivial length 1 composition) 
and set
${\cal R}:=\bigcup_{m \in \Z} {\cal R}_m$.

We write $res$ for the restriction homomorphism $Aut(V) \rightarrow GL(U)$.  Since
$U$ generates $V$, $res$ is injective.  We consider the question of when $g \in
GL(U)$ is in the image of $res$.

We shall define a set function $e:GL(U) \rightarrow End(V)$ as
follows. For $g \in GL(U)$, define $e(g) \in End(V)$ by its
action on the basis elements $\mu (\vec x)$, $(\mu,\vec x) \in {\cal R}:$
\begin{equation}\label{5}
e(g)(\mu (\vec x)):= \mu (g(\vec x)). 
\end{equation}
This endomorphism will turn out to be invertible in cases of interest to us.

Now consider the following set of conditions on $e(g) \in End(V)$.

\begin{equation}\label{7.a}
e(g)(\mu (\vec u)) = \mu (g(\vec u))
\end{equation}
\begin{equation}\label{7.b}
e(g)e(g^{-1})(\mu (\vec u)) = \mu (\vec u) = e(g^{-1})e(g)(\mu (\vec u)), 
\end{equation} 
for all $U$-admissible pairs $(\mu,\vec u).$

We may assume that the components of $\vec x$ are homogeneous elements,
and even $(\mu, \vec u) \in {\cal Q}_m$.  Both sides of 
(\ref{7.a}) are expanded 
in the basis ${\cal B}_m$.   
Equating the coefficients of both sides gives {\it polynomial
conditions} on the entries of $(g_{ij})$, the  matrix representing $g$ with respect
to $\cal B$. A similar discussion applies to (\ref{7.b}).

There is an ideal $I_{(\mu,\vec x)}$ in the 
ring 
$\C [ x_{ij}, det^{-1}  
|  i,j = 1, \dots , dim(U)  ]$ 
of polynomial functions on $GL(U)$ associated to the above 
set of conditions. 

Finally, for $u \in U$, define the ideal $I_u$ by the condition $g u = u$.  
Set $I :=\sum _{(\mu,\vec x) \in \cal Q} I_{(\mu,\vec x)} + I_{\vac}
+I_\omega$ and
$$G_U= \{ g \in GL(U) \ | \ p(g)=0, \hbox { for all } p \in I \}.$$
Then $G_U$ is a variety contained in $GL(U)$.   
Clearly, $res(Aut(V)) \le G_U$.  

\begin{lem} $G_U$  is a subgroup of $GL(U)$ defined 
as the zero set for the ideal $I$.   
That is, $G_U$ is an algebraic group.  Also, $e$ is a homomorphism.  
\end{lem}

\pf  First, $1 \in G_U$.  Observe that for $g \in G_U$, $e(g)$  and $e(g^{-1})$ are
invertible since their restrictions to each $V_i$ are invertible. 
Moreover, they form
an inverse pair, whence $e(g^{-1})=e(g)^{-1}$. 
We now show that $g^{-1}$ satisfies (\ref{7.a}). 
Let $(\mu, \vec y)$ be a $U$-admissible pair of length $r+1$.  
Since $g \in G$, we
have, for all $\vec y$, $e(g)\mu (g^{-1}(\vec y) ) = \mu (\vec y)$ and so 
$\mu (g^{-1}(\vec y )) = e(g)^{-1} \mu (\vec y) = e(g^{-1}) \mu (\vec y)$. 
Since
$g^{-1}$ satisfies (\ref{7.a}), we have $G_U=G_U^{-1}$. 

To prove closure under products, we let 
$(\mu, \vec u)$ be a $U$-admissible pair and $g, h \in
G_U$.  We must show that 
$e(gh)\mu (\vec u)=\mu (gh(\vec u))$.   Write $\mu (\vec u)=\sum_{(\nu,\vec y) \in \cal R} a_{(\nu,\vec y)}
\nu(\vec y)$, for unique scalars $a_{(\nu,\vec y)}$, almost all zero. 
Then,
$$e(gh)\mu (\vec u)=\sum_{(\nu,\vec y) \in \cal R} a_{(\nu,\vec y)}
e(gh)\nu(\vec y) =   
\sum_{(\nu,\vec y) \in \cal R} a_{(\nu,\vec y)} \nu(gh(\vec y)),$$  
by definition of $e(gh)$.  Since $g \in G_U$, this equals
$$\sum_{(\nu,\vec y) \in \cal R} a_{(\nu,\vec y)} e(g) \nu(h(\vec y)) = 
\sum_{(\nu,\vec y) \in \cal R} a_{(\nu,\vec y)} e(g)e(h) \nu(\vec y) = e(g)e(h)\mu (\vec u).$$
Also $\mu (gh(\vec u))=e(g)\mu (h(\vec u))=e(g)e(h)\mu (\vec u)$ because $g, h \in G_U$.  
We conclude that $gh \in G_U$
and so $G_U$ is a group.  

	Since the $\mu (\vec u)$ span $V$, we 
also deduce that
$e(gh)=e(g)e(h)$, whence $e$ is a homomorphism.    
\qed

\begin{lem} For all $u, v \in V$, $n \in \Z$, we have
$$e(g)(u_nv)=(e(g)u)_n(e(g)v).$$
That is, $Im(e) \subseteq Aut(V)$.  
\end{lem}

\pf  We may assume that $u$ is ``monomial", i.e., has the form 
$\mu (\vec x)$,
for a $U$-admissible pair $(\mu, \vec x)$.   We argue by induction on the
length of 
$(\mu,\vec x)$. First, we assume that the  length is 1.  We may also
assume  that $v$ is monomial, so $v=\nu(\vec y)$, for a $U$-admissible
pair $(\nu , \vec y)$.  Say $\nu$ is an $\vec m$-composition, $\vec m =
(p , \dots, q)$  and
$\vec y = (y_1, \dots , y_t)$.  Then,  
$u_nv= u_n (y_1)_p\cdots {}_q(y_t)$ and 
$e(g)(u_nv)=  (gu)_n (gy_1)_p \cdots {}_q(gy_t)$, by \ref{7.a}.  
By \ref{7.a} applied to $(\nu , \vec y)$, we deduce 
$e(g)(u_nv)=  (gu)_n (e(g)\nu (\vec y)) = (gu)_n (e(g)v)$.  
Finally, since $e(g)x=gx$ for $x \in U$, this is $(e(g)u)_n(e(g)v)$.

Suppose next that the length is $r\ge 2$ and that $\mu$ 
is an $\vec m$-composition, $\vec m = (m_1, \dots , m_r)$. 
Set $k=m_1, b=x_1$,
$a=\nu(\vec y)$, where $y=(x_2, \dots , x_r)$ and 
$\nu$ is the  $V$-composition
associated to the $(r-1)$-tuple $(m_2, \dots , m_r)$. Then $u=b_ka$.

We now do a residue calculation to verify that 
$$e(g)(Y(u,z)v)= Y(e(g)u, z)(e(g)v).$$
Extracting the coefficient at $z^{-n-1}$ will give the Lemma.

Since $u=b_ka$, we have from the Jacobi identity for vertex operators
(see the formula before (3.3) of [D]) that 
$$Y(u,z)v = \Res_w\{ (w-z)^k
Y(b,w)Y(a,z)v - (-z+w)^kY(a,w)Y(b,z)v \}.  $$

Write $h$ for $e(g)$. Then 
$$h[ Y(u,z)v ]\!=\!\Res_w \{ (w-z)^k h[Y(b,w)Y(a,z)v] -
(-z+w)^kh[Y(a,w)Y(b,z)v]\}.$$ 
By several uses of 
induction on length, applied to $b$ and $a$, plus the above consequence
of the Jacobi identity, 
we deduce that 
this equals
\begin{eqnarray*}
 & &\ \ \ \Res_w\{ (w-z)^k [Y(hb,w)h[Y(a,z)v] -
(-z+w)^k[Y(ha,w)h[Y(b,z)v]] \}\\
& & =\Res_w\{ (w-z)^k [Y(hb,w)Y(ha,z)](hv) \\
& &\ \ \ \ -
(-z+w)^k[Y(ha,w)Y(hb,z)(hv)] \}\\
& & =Y((hb)_k(ha),z)(hv)\\
 & & =Y(h(b_ka),z)(hv)\\
 & &=Y(hu,z)(hv),
\end{eqnarray*}
as desired.   
\qed

\begin{thm}
The two maps 
$$res: Aut(V)
\rightarrow G_U \ \ \hbox{ and } \ \ e:G_U \rightarrow Aut(V)$$ 
form a pair of inverse
isomorphisms. 
Therefore, $Aut(V)$ is isomorphic to the algebraic group $G_U$.
\end{thm}

\pf Since $U$ generates $V$, $res$ is a monomorphism. Since $Im(e)$ is
contained in $Aut(V)$ and $res \circ e = Id_{G_U}$,  $res$ is an epimorphism hence an
isomorphism.   Since the set map $e$ is a one sided inverse of an isomorphism (hence
a two-sided inverse),  it is an isomorphism of groups.  (We proved before
that $e$ is a homomorphism, but do not need to quote that result here.)
\qed

\begin{rem}{\rm The most well known vertex operator algebras are finitely
generated. For examples, Heisenberg vertex operator algebras [FLM] and 
affine vertex operator algebras (cf. [DL], [FZ], [L2]) are generated by 
their
weight one subspaces, Virasoro vertex operator algebras (cf. [FZ], [L2]) 
and the moonshine vertex operator algebra (see [B], [FLM]) are 
generated by weight 2 subspaces.  The lattice vertex operator
algebra
$V_L$ (see [B], [FLM]) is generated by 
$\oplus_{m\leq n}(V_L)_m$ where
$n$ is any positive integer such that 
$L$ has a direct sum decomposition 
$L=\oplus_{i=1}^n \Z \alpha_i$ satisfying  
$\la \alpha_i, \alpha_i \ra /2 \leq n$. In fact, $V_L$ is generated by
$e^{\pm \alpha_i}$ for $i=1,...,n.$  }
\end{rem} 

\begin{ex}\label{exama}{\rm 
If $G$ is not finitely generated, $Aut(V)$ is not an algebraic group
in general. Here is an example. Let $(U,Y,{\bf 1},\omega)$ be a
vertex operator algebra with infinitely many irreducible modules
$U^i=(U^i,Y^i)$ ($i=1,.2,...$) not isomorphic to $U$ such that
$U^i=\oplus_{n\geq 0}U^i_{\l_i+n}$ with $U^i_{\l_i}\ne 0$ and
$\l_1<\l_2<\cdots.$ Set $V=U\oplus \bigoplus_{i>0}U_i.$ Then $V$  
has a vertex
operator algebra structure 
with vertex operator $Y'$ defined in the following way (see [L1]).
Since 
$V$ is a $U$-module, $Y'(u,z)v$ for $u\in U$ and $v\in V$ is defined 
in an obvious way. Using the idea of skew symmetry, for $v \in U^i$ and 
$v \in U$,   we define 
$Y'(v,z)u:=e^{zL(-1)}Y'(u,-z)v.$
Finally we define $Y'(v,z)w:=0$ for all $v,w\in \oplus_{i}U^i.$ We refer
the reader to [L1] for the proof that $(V,Y',{\bf 1}, \omega)$ 
is indeed a vertex operator algebra.  

For $k,\lambda\in \C$ let $L(k,\lambda)$ be the irreducible highest weight
module for the Virasoro algebra with central charge $k$ and
highest weight $\l.$ Then $L(1,0)$ is a vertex operator 
algebra and $L(1,\lambda)$ is an irreducible $L(1,0)$-module
for any $\l$ (cf. [FZ]). Now we take $U=L(1,0)$ and $U^i=L(1,i)$, 
for $i = 1, 2, \dots $.   
Let $u^i$ be a nonzero highest weight vector of $U^i$ (which is unique
up to a scalar). Then $V$ is generated by $\omega$ and $u^i$ for
$i>0.$ Clearly, $V$ is not finitely generated since a finite set of generators
would lie in the sum of $U$ and finitely many $U^i$.  Note also that the
sum of any set of the $U^i$ is an ideal.   }
\end{ex}

\begin{prop} The automorphism group of the VOA $\bigoplus_{n=0}^\infty
L(1,n)$ is isomorphic to the infinite direct product
$\prod_{i=1}^{\infty}\C^{\times}_i$ where
$\C^{\times}_i$ is a copy of multiplicative group $\C^\times$ acting 
faithfully on $U^i$, trivially on $U^j$ for $j \ne i$ and trivially
on $U$. In particular, $Aut(V)$ is not an algebraic group. 
\end{prop}

\pf Let $\l=(\l_1,\l_2,...)\in  \prod_{i=1}^{\infty}\C^{\times}   $.  
We define a $U$-module homomorphism
$g_{\l} \in \prod_{i=1}^{\infty}\C^{\times}_i$ on 
$V := \bigoplus_{n=0}^\infty L(1,n)$ 
such that
$g_{\l}{\bf 1}={\bf 1}$ and
$g_{\l}u^i=\l_iu^i.$ It is easy to see from the definition of $Y$ that
$g_{\l}$ is an automorphism of $V.$ On the other hand any automorphism
$g$ is identity on $U$ as $U$ is generated by the Virasoro element. So
$g$  preserves the space of highest weight vectors which is spanned by
${\bf 1}$ and $u^i$ for $i>0.$ Since the weights of any two 
highest weight vectors are different, we immediately have that $gu^i=\l_iu^i$
for some nonzero constant $\l_i$ for all $i.$ Set $\l=(\l_1,\l_2,\cdots)
\in \prod_{i=1}^{\infty}\C^{\times}.$ Then $g=g_{\l}.$ 
Clearly, $\l \mapsto g_\l$ is an isomorphism.  \qed

Next we discuss the automorphism group of $V$ for a ``nice''
vertex operator algebra. We need more definitions. 

\begin{de} {\rm A vertex operator algebra
 $V$ has {\em CFT type} if $V_n=0$ for $n<0$ and
$dim(V_0)=1$ (so $V_0=\C \vac$).}
\end{de}

In the following definition we use the notion of admissible modules
as introduced in [Z] and [DLM2]. We refer the reader to [DLM2]
for details.

\begin{de}{\rm A vertex operator algebra $V$ is {\it rational} 
 if any admissible
module is a direct sum of irreducible admissible modules.}
\end{de}

\begin{de}{\rm A vertex operator algebra $V$ is $C_k$-{\it cofinite} 
if $\dim (V/C_k(V))$ 
is finite where $C_k(V)$ is the  subspace
of $V$ spanned by $u_{-k}v$ for $u,v\in V.$ }
\end{de}

The $C_2$-cofinite condition in the literature was called the 
$C_2$-finite
condition or $C_2$-condition (as in [Z]). In the case of vertex operator 
algebras associated to highest weight modules for affine Lie algebras
and the Virasoro algebras, $V/C_k(V)$ are the spaces of  
coinvariants (cf. [FF], [FKLMM]).  This should explain why we are 
changing the terminology. In this paper only the $C_2$-cofinite condition is 
used.

As we have already mentioned that if $V$ is of CFT type then $(V_1,
0^{th})$ is a Lie  algebra under $[u,v]=u_0v.$ 
The part (i) of the 
following theorem can be found in [DM2].
The rest follows from the general structure of algebraic groups.
\begin{thm}\label{tdm}
Let $V$ be a simple, $C_2$-cofinite rational  vertex operator algebra of CFT 
type with $L(1)V_1=0.$ Then,
\begin{enumerate}
\item  $V_1$ is a reductive Lie algebra; 
write $V_1=  {\mathfrak s} \oplus {\mathfrak t}$, where the first summand 
is semisimple and the second is toral.
\item	$G:=Aut(V)$ contains the connected component $G^0$ of the
identity with finite index and satisfies 
$G^0=G_1C_1$ (central product) where $G_1 := 
\la exp(x_0) \ | \ x \in V_1 \ra$ and  $C_1 := C_G(V_1)^0$.
We have $G_1=S_1T_1$, where  
$S_1 := \la exp(x_0) \ | \ x \in {\mathfrak s} \ra$  
and 
$T_1 := \la exp(x_0) \ | \ x \in  {\mathfrak t} \ra$.  
Also, $T_1 = (C_1 \cap G_1)^0$ and 
there is a connected group $K_1$ which is normal in $G$
and has the properties
 $C_1=T_1K_1$, $[T_1,K_1]=1$ and $T_1 \cap K_1$ is finite.  
\end{enumerate}  
\end{thm}

We remark that the condition $L(1)V_1=0$ is not a strong assumption. 
It seems that all rational vertex operator algebras
of CFT type satisfy this condition. For example, 
well known rational vertex operator algebras associated to
highest weight integral modules for affine algebras (cf. [L2]),
 to minimal series for the Virasoro algebras (the weight one space
is zero in this case), or to  positive definite even lattices
satisfy the condition (cf. [FLM]). It is proved in [L1] that for a 
simple vertex operator algebra of CFT type,
the condition $L(1)V_1=0$ is equivalent to that there is a 
nondegenerate symmetric invariant bilinear form on $V$ in the sense
of [FHL].  

\section{Derivations}

There is a close relation between automorphisms and derivations
for a vertex operator algebra. In this section we discuss the Lie algebra
of the derivations of a vertex operator algebra. 

Define a linear map $o$ on $V$ by setting $o(v)=v_{\wt v-1}$ for homogeneous elements $v$. 
Then $o(v)V_n\subset V_n$ for all $n$.  

A {\it derivation } of the vertex operator algebra $V$ is an endomorphism $d$ of
$V$ which satisfies $d\vac = 0, d\o = 0$ and $[d,Y(u,z)]=Y(du,z)$. 
Since $d\o=0$, $d$ preserves all the $V_n$, which are the eigenspaces
of an operator in $Y(\o ,z)$, whence $d$ is locally finite.  
The derivation $d$
is an {\it inner derivation} if there is $v \in V$ so that $o(v) = d$ 
(see \ref{oV1}).  
Since the $V_k$ are finite dimensional, any endomorphism preserving
the graded pieces is locally finite.  

Since $d$ is a locally finite derivation of $V$, the exponential
$e^d$ is an automorphism of $V.$ On the other hand, $Aut(V)^0$
(when $V$ is finitely generated) 
is a connected Lie group and its Lie algebra acts on $V$ as derivations.

Set $IDer(V):= o(V)\cap Der(V)$, the space of {\it inner derivations}.  

Let $V$ be of CFT type such that $L(1)V_1=0$. 
Then
$V$ is a direct sum of irreducible modules for $span\{ L(-1),L(0), L(1) \}
\cong sl(2,\C)$,  the {\it
principal $sl_2$} [DLinM]. For homogeneous $v$, 
since $o(L(-1)v)= -(wt(v)-1)v_{wt(v)-1}$, we have
equality of $\{ o(v) \ | \ v \in V \}$ and $\{ o(v) \ | \ v \in Ker(L(1))\ \}$.

Let $QV : = Ker(L(1))$, the space of {\it quasi primary vectors}.

\begin{lem}\label{oV1}  We have 
$o(v) \in Der(V)$ for $v \in V_1$.
\end{lem}

\pf Since $[o(v),Y(u,z)]=[v_0,Y(u,z)]=Y(v_0u, z)=Y(o(v)u,z)$ for
$v\in V_1$ and $u\in V,$ the result is clear.
\qed

\begin{lem} Assume that $V$ has CFT type. If $v = \sum_{i \ge 2} v^i$, 
with $v^i \in V_i\cap QV$ and $o(v)
\in Der(V)$, then $v=0$.
\end{lem}

\pf  Since $o(v)=\sum_{i \ge 2} v_{i-1}^i$, we have 
\begin{eqnarray*} 
& &[o(v),Y(u,z)] =\sum_{i \ge 2}[v_{i-1}^i, Y(u,z)]\\
& &\ \ \ \ \  = \sum_{i \ge 2}\sum_{j \ge 0} \binom{i-1}{j}
Y(v_j^iu,z)z^{i-1-j}\\
& &\ \ \ \ = \sum_{i \ge 2}Y(v_{i-1}^iu,z).
\end{eqnarray*}
It follows that 
$$\sum_{i \ge 2}\binom{i-1}{j}\sum_{j=0}^{i-2}Y(v_j^iu,z)z^{i-1-j}=0$$
 and 
$$lim_{z \rightarrow 0}\{ \sum_{i \ge
2}\binom{i-1}{j}\sum_{j=0}^{i-2} Y(v_j^iu,z)z^{i-1-j-1}\vac \}  = 0.$$

This implies 
$$\sum_{i \ge 2} \binom{i-1}{i-2}v_{i-2}^i u = 0,$$ for all $u.$
Thus 
$$\sum_{i \ge 2}\binom{i-1}{i-2}v_{i-2}^i=0$$ on $V.$ Since the 
$v^i$ are quasi-primary vectors, 
$[L(1),v_{i-2}^i]=iv_{i-1}^i$. As a result
we have
$$\sum_{i \ge 2}(i-1)iv_{i-1}^i=0.$$

By Theorem 2.2 of [DLMM], we get $\sum_{i \ge 2}i(i-1)v^i \in V_1$, whence $v^i=0$
for all $i$ and so $v=0$.
\qed

The following corollary is immediate:
\begin{coro}  For $V$ of CFT type, 
$IDer(V) = o(V_1)=\{o(v)|v\in V_1\}.$
\end{coro}

Recall from [DLMM] that  the radical $J(V)$ of $V$ consists 
of those vectors $v\in V$ such that $o(v)=0.$ 
We need a result from [DM2].
\begin{lem}\label{la1} Let $V$ be a $C_2$-cofinite 
rational vertex operator algebra of CFT 
type. Then $J(V)=(L(-1)+L(0))V.$
\end{lem}

{}From now on we assume that $V$ is a $C_2$-cofinite 
rational vertex operator algebra of CFT type. 
Then $\mathfrak g=(V_1, 0^{th})$ is a reductive 
Lie algebra and each $V_n$ is a finite dimensional $\mathfrak g$-module
via $v\mapsto o(v).$ Define invariant symmetric bilinear form 
$(,)_M$ on $\mathfrak g$ for any $\mathfrak g$-module $M$ by 
$(u,v)_M=\tr_{M}(uv)$ for $u,v\in \mathfrak g.$ 

Recall that $\mathfrak g=\mathfrak s\oplus \mathfrak t$ where $\mathfrak s$ is 
semisimple and $\mathfrak t$ is abelian. Then each finite dimensional module
for $\mathfrak g$ is a direct sum of indecomposable modules
which are tensor products of irreducible modules for 
${\mathfrak s}$ and indecomposable modules for $\mathfrak t.$ 

\begin{lem}\label{la2}
 Let $M$ be a finite dimensional $\mathfrak g$-module such that
$M$ contains $\mathfrak s$ as an $\mathfrak s$-module. Let ${\mathfrak s}={\mathfrak s}^1\oplus 
\cdots \oplus {\mathfrak s}^p$ be the decomposition into simple ideals. 
Write ${\mathfrak s}^0:= {\mathfrak t}$.  
Then
$({\mathfrak s}^i,{\mathfrak s}^j)_M=0$ if $i\ne j$ and if $i > 0$, 
the restriction of the form to each ${\mathfrak s}^i$ is nondegenerate.
\end{lem}

\pf First we prove that the restriction of  the form 
to each ${\mathfrak s}^i$ is nondegenerate. Note that as an ${\mathfrak s}^i$-module,
$M$ is completely reducible and ${\mathfrak s}^i$ is an irreducible submodule.

Let $i > 0$.  
It is well known that for each irreducible ${\mathfrak s}^i$-module $W$
the corresponding invariant symmetric bilinear form $(,)_W$ is
a nonnegative multiple of $(,)_{{\mathfrak s}^i}$, the Killing form
on ${\mathfrak s}^i$, which is nondegenerate (cf. [H]). As a result, $(,)_M$ 
is nondegenerate when restricted to ${\mathfrak s}^i.$

In order to prove that $({\mathfrak s}^i,{\mathfrak s}^j)_M=0$ if $i\ne j$, 
we may assume that $M$ is irreducible.   
Then $\mathfrak t$ acts as scalars on $M$ and 
$M|_{\mathfrak s} =M^0\otimes \cdots  \otimes M^p$ is a tensor product of
irreducible modules $M^i$ for ${\mathfrak s}^i$.  

Let $i \ge 0$ be any index and let  $x\in {\mathfrak s}^i.$ 
Then $M^i$ is a direct sum of generalized 
eigenspaces under $x.$ So, we can chose a basis $B_k$ for $M^k$ 
consisting of generalized eigenvectors for 
the action of $x$ on  $M^i.$ Let $m^k \in B_k$, 
associated to the  
generalized eigenvalue $\lambda_k.$  
Let $\lambda := \sum_{k \ne j} \lambda _k$.   
Now fix $j >0$ and assume 
$j \ne i$.  Let $y\in {\mathfrak s}^j$.   
Think of a  matrix for the action of $y$ written
in a basis taken from 
the subspaces of the form 
$m^0\otimes m^1 \otimes  \cdots \otimes 
m^{j-1}  \otimes M^j \otimes m^{j+1}  
\cdots \otimes m^p$.   
Note that $M$ can be written as a direct sum of such subspaces.  
These spaces are not invariant by $x$ but the action  of $xy$ has
block-triangular form with respect to this direct sum.  

The contribution to the trace $\tr_{M}xy$ from 
the subspace 
$m^0\otimes m^1 \otimes  \cdots \otimes 
m^{j-1}  \otimes M^j \otimes m^{j+1}  
\cdots \otimes m^p$
is equal to
$\lambda\tr_{M^j}y$ for $y\in {\mathfrak s}^j.$ Since ${\mathfrak s}^j$ is
simple and $[{\mathfrak s}^j,{\mathfrak s}^j]={\mathfrak s}^j$ we see that
$\tr_{M^j}y=0.$ Thus  $({\mathfrak s}^i,{\mathfrak s}^j)_M=0$ since $M$
is a direct sum as described above.  
\qed

For convenience we denote the bilinear form $(,)_{V_n}$ on
$\mathfrak g$ by $(,)_n$ for $n\geq 0.$ 

We need the following result from [DM2].

\begin{prop}\label{impo} Let $V$ be a rational, $C_2$-cofinite
vertex operator algebra of CFT type and $L(1)V_1=0.$ Then for
any $u\in V_1$ there exist $n\geq 0$ and $v\in V_1$ such that
$(u,v)_n\ne 0.$
\end{prop}

The next result sharpens Proposition \ref{impo}.

\begin{thm}\label{la3} Let $V$ be as in Proposition \ref{impo}. Then there exists 
$n$ such that $(,)_n$ is nondegenerate.
\end{thm}

\pf Take $n$ large enough so that $\sum_{m=0}^nV_m$ generates $V.$ 
We claim that $(,)_n$ is nondegenerate. 

Recall that $\C L(-1)+\C L(0)+\C L(1)$ is isomorphic to Lie algebra
$sl(2,\C).$ Let $M(\lambda)$ is the irreducible
highest weight module for $sl(2,\C)$
with highest weight $\lambda.$ Then 
$$M(\lambda)=\oplus_{m\geq 0}M(\lambda)_{\lambda+m}$$
such that $M(\lambda)_{\lambda}$ is spanned by a highest weight vector
$v_{\lambda}$ and $M(\lambda)_{\lambda+m}$ is spanned by
$L(-1)^mv_{\lambda}.$ If $\lambda=0$ then $M(\lambda)$ is trivial.
If $\lambda>0$ each $M(\lambda)$ is a Verma module and
$M(\lambda)_{\lambda+m}\ne 0$ for all $m.$ 

First we prove that the representation of $\mathfrak g$ on $V_n$ is faithful.
Assume that $u_0=0$ on $V_n.$   Since $V$ is CFT type and
$L(1)V_1=0$ it follows from Corollary 3.2 of [DLinM] that  
$V$
is a direct sum of copies of 
 $M(\lambda)$ for $m\geq 0$ 
 and the multiplicity of
$M(0)$ in the decomposition is 1. Note that $[L(i),u_0]=0$ for
$i=-1,0,1$ and $u\in \mathfrak g.$ Let $M(\lambda)$ occurs in the decomposition
of $V$ such that $0\ne m\leq n.$ Then $u_0$ on $M(\lambda).$ Also
note that $u_0V_0=0.$ Thus $u_0=0$ on $\oplus_{m=0}^nV_m.$ Since $u_0$
is a derivation on $V$ and $V$ is generated by $\oplus_{m=0}^nV_m,$ we immediately see that $u_0=0$ on $V.$ This contradicts Proposition \ref{impo}.

So we can identify $\mathfrak g$ with its image $o(\mathfrak g)_n=
\{u_0|_{V_n}\}.$ If the form $(,)_n$ is degenerate then by Lemma
\ref{la2}, there exists $x\in {\mathfrak t}$ such that $(x,y)_n=0$ for
all $y\in \mathfrak g.$ By Lemma 4.3 of [H], $x_0$ is nilpotent on $V_n.$
In particular, all the eigenvalues of $x_0$ on $V_n$ are zero. 
A similar argument as above then shows that $x_0$ has only zero 
eigenvalues on $\oplus_{m=0}^nV_m.$ Since $\oplus_{m=0}^nV_m$  generates $V$ we
see immediately that $x_0$ has only zero eigenvalues
on $V_m$ for all $m$
as 
$$u_0v^1_{m_1}\cdots v^k_{m_k}v=\sum_{j=1}^kv^1_{m_1}\cdots (u_0v^j)_{m_j}\cdots v^k_{m_k}v+v^1_{m_1}\cdots v^k_{m_k}u_0v$$
for $v^j,v\in V$ and $m_j\in \Z.$ 

Note that ${\mathfrak t}$ is an abelian Lie algebra and all irreducible
modules are one-dimensional. So,  $\tr_{V_m}x_0y_0=0$ for
all $y\in \mathfrak g$ and $m\in \Z.$ Again, by Proposition \ref{impo}, this
is impossible.  \qed

\begin{thm}
Let $V$ be a $C_2$-cofinite rational vertex operator
algebra of CFT type such that $L(1)V_1=0$.  
Let $n\geq 0$ such that $\sum_{m=0}^n V_m$ generates $V.$
Then $Der(V)$ is a direct sum of ideals 
$o(\mathfrak g)$ and ${\mathfrak g}^{\perp}$
where ${\mathfrak g}^{\perp}$ consists of $d\in D$ such that
$\tr_{V_n}o(u)d=0$ for all $u\in V_1.$ 
\end{thm}

\pf  
Let $V$ be as in Lemma \ref{la2} and Theorem \ref{la3}. Let $n>0$ be as in the
proof of Theorem \ref{la3}. Then the action of $D:=Der(V)$ is also faithful
on $V_n$ and $(d,d')_n=\tr_{V_n}d d'$ defines a symmetric 
invariant bilinear form on $D.$ 
Then restriction of the form to
$o(\mathfrak g)$ 
is nondegenerate by Theorem \ref{la3}. Let $\mathfrak g^{\perp}$ be
the orthogonal complement of $o(\mathfrak g).$ Then the intersection
of ${\mathfrak g}^{\perp}$ and $o(\mathfrak g)$ must be zero. On the other
 hand,
$[d,u_0]=(du)_0$ tells us that $o(\mathfrak g)$ is an ideal of $D$
and so is ${\mathfrak g}^{\perp}.$ Thus $[d,u_0]=(du)_0=0$
for $d\in {\mathfrak g}^{\perp}$ and $u\in V_1.$ Since the action of
$D$ is faithful on $V_n$ we have $du=0.$
\qed

\section{Example}

In this section we give an example 
which shows that a finitely generated 
VOA with an infinite descending chain of 
ideals can still have a reductive group of automorphisms.  
Our example is $V_L^G$, below.  
In this example, we find all ideals and find that they form
a countable descending chain.

\begin{ex} \rm  
We consider the vertex operator algebra $V=V_L=L(\Lambda_0)$
where $L=\Z\alpha$ such that $(\alpha,\alpha)=2$ and $L(\Lambda_0)$
is the fundamental representation for the affine algebra $A_1^{(1)}.$
Then the automorphism group of $V_L$ is isomorphic to
$PSL(2,\C)$ (see [DLY] and [DN]). 

Let $L(c,h)$ be
the highest weight module for the Virasoro algebra with central charge
$c$ and highest weight $h.$ Let $W_m$ for $m\geq 1$ be
the irreducible module for $sl(2,\C)$ of dimension $m.$ Then
$$V_L=\oplus_{m\geq 0}L(1,m^2)\otimes W_{2m+1}$$
(cf. [DG]) and $SL(2,\C)$ acts on $V_L$ by acting on $W_{2m+1}.$
Moreover, $W_{2m+1}$ regarded as a $SL(2,\C)$-submodule
of $V_L$ is generated by the highest weight vector $e^{m\alpha}.$

Consider the subgroup 
$G=\left\{\left(\begin{array}{cc} 1 
& t\\ 0 & 1\end{array}\right)|t\in \C \right\}$ 
of $SL(2,\C).$  Clearly $G$ is not compact.
Then the space of $G$-invariants $W_{2m+1}^G$ is spanned by
$e^{m\alpha}.$ 
As a result, 
we have a direct sum decomposion of the fixed point set for $G$ into 
irreducible modules for the Virasoro algebra: 
$$V_L^G=\oplus_{m\geq 0}L(1,m^2),$$ where the  
highest weight module
$L(1,m^2)$ for the Virasoro algebra is generated by $e^{m\alpha}.$ 
For distinct $m$, these are pairwise nonisomorphic. 
It is
not hard to see that for any $n\geq 1,$ $\sum_{m\geq n}L(1,m^2)$ is an
ideal of $V_L^G$ as $u_kv\in L(1,(s+t)^2)$ for any $u\in L(1,s^2)$ and
$v\in L(1,t^2)$ and $k\in\Z.$   We now prove that 
all the ideals of
$V_L^G$ are
given in this way. Let $I$ be a nonzero ideal of $V_L^G.$ Then $I$ is
a module for the Virasoro algebra, so is a sum of a family of the
above $L(1,n^2)$. 
Let
$n\geq 1$ be the smallest
positive integer such that $L(1,n^2)$ is a subspace of $I.$ Then
$e^{\alpha}_{-n(\alpha,\alpha)-1}e^{n\alpha}=e^{(n+1)\alpha}\in I$
(see [FLM]) and $L(1,(n+1)^2)$ is contained in $I.$ 
Continuing in this way we see that
$I=\sum_{m\geq n}L(1,m^2).$ 
It was proved in [DLM1] that if $H$
is a compact Lie group which acts continuously on a simple vertex
operator algebra $V$ then $V^H$ is also a simple vertex operator
algebra.  Since our $G$ is not compact, $V_L^G$ is permitted to be
nonsimple, and it is.

We now show that the automorphism group of $V_L^G$ is isomorphic
to $\C^\times.$ Note that $V_L^G$ is generated by $\omega$ and $e^{\alpha}.$
In fact, we have already proved that $e^{(n+1)\alpha}$ can be generated
from $e^{\alpha}$ and $e^{n\alpha}.$ 
So all the highest weight vectors
can be generated from $e^{\alpha}.$ Using $\omega$ then generates the
whole space. For any $\lambda\in \C^\times$ we define
a linear isomorphism $\sigma_{\lambda}$ of $V_L^G$ such that
$\sigma_{\lambda}$ acts on $L(1,m^2)$ as $\l^m.$ From the discussion
before it is clear that $\sigma_{\l}$ is an automorphism
of $V_L^G.$ On the other hand, any automorphism $\sigma$ of $V_L^G$
maps $e^{\alpha}$ to $\lambda e^{\alpha}$ for some
nonzero constant $\l$ as $e^{\alpha}$ is the only highest weight
vector with highest weight 1 (up to a constant)
for the Virasoro algebra. Since $V_L^G$ is generated by $\omega$ and
$e^{\alpha}$ we immediately see that $\sigma=\sigma_{\lambda}.$
\end{ex}

\end{document}